\input ./style/arxiv-general.cfg
\documentclass[aop,MSNbibl,dvips]{arximspdf}
\makeatletter
   \@ifpackageloaded{graphicx}{}{\usepackage{graphicx}}
\makeatother

\doi{10.1214/14-AOP982}
\volume{44}
\issue{1}
\pubyear{2016}
\firstpage{567}
\lastpage{588}
\docsubty{FLA}

\makeatletter
\newtheorem{hypo}{Hypothesis}

\newcommand{\rrvert}{\vert}

\newcommand{\llvert}{\vert}
\newtheorem{theorem}{Theorem}[section]
\newtheorem{proposition}[theorem]{Proposition}
\newtheorem{corollary}[theorem]{Corollary}
\newproclaim{remark}[theorem]{Remark}
\newproclaim{definition}[theorem]{Definition}
\newtheorem{lemma}[theorem]{Lemma}

\newcommand{\R}{\mathbf R}
\newcommand{\E}{\mathbf E}

\makeatother

\begin{document}
\begin{frontmatter}

\title{An averaging principle for diffusions in~foliated~spaces}
\runtitle{Averaging principle in foliated spaces}

\begin{aug}
\author[A]{\fnms{Ivan I.}~\snm{Gonzales-Gargate}\thanksref{T1}\ead[label=e1]{ivangargate@utfpr.edu.br}}
\and
\author[A]{\fnms{Paulo R.}~\snm{Ruffino}\corref{}\thanksref{T2}\ead[label=e2]{ruffino@ime.unicamp.br}}
\runauthor{I.~I. Gonzales-Gargate and P.~R. Ruffino}
\thankstext{T1}{Supported by CNPq 555241/2009-2.}
\thankstext{T2}{Supported in part by FAPESP 11/50151-0 and 12/03992-1.}
\affiliation{Universidade Estadual de Campinas}
\address[A]{Departamento de Matem\'{a}tica\\
Universidade Estadual de Campinas\\
13.083-859- Campinas - SP\\
Brazil\\
\printead{e1}\\
\phantom{E-mail:\ }\printead*{e2}}
\end{aug}

%
\received{\smonth{1} \syear{2013}}
%
\revised{\smonth{3} \syear{2014}}

%
\begin{abstract}
Consider an SDE on a foliated manifold whose trajectories lay on
compact leaves. We investigate the effective behavior of a small
transversal perturbation of order $\varepsilon$. An average principle is
shown to
hold such that the component transversal to the leaves converges to
the solution of a deterministic ODE, according to the average of the
perturbing vector field with respect to invariant measures on the
leaves, as
$\varepsilon$ goes to zero. An estimate of the rate of convergence is
given. These
results generalize the geometrical scope of previous approaches, including
completely integrable stochastic Hamiltonian system.
\end{abstract}

%
\end{frontmatter}

\section{Introduction and set-up}\label{sec1}

Generally speaking, the original heuristic idea of an averaging principle
refers to an intertwining
of two dynamics where one
of them is, in some sense, much slower and is affected
somehow by the other faster dynamics. An
averaging principle in this case
refers to the possibility of approximating, in some topology, the slow dynamics
considering only the
average action or perturbation which the fast motion induces on it.
These ideas appeared long ago, and as mentioned by Arnold (\cite{V-Arnold},
page
287),
they
were implicitly contained in the works of Laplace, Lagrange and Gauss on
celestial mechanics; literature on the matter can be found, for example,
\cite{V-Arnold}, Sanders,
Verhulst and
Murdoch \cite{SVM} and references therein, among many others.
Presently, regarding
stochastic systems, averaging has been quite an active
research field on which there is also a vast literature on the topic.
Interesting quick historical overviews can
be found in Li (\cite{Li}, page 806), Kabanov and Pergamenshchikov
(\cite{Kabanov-Pergamenshchikov}, Appendix), (\cite{SVM}, Appendix~A).
Among many
other works somehow related to the topic, we refer to
Khasminski and
Krylov \cite{Khasminski-krylov}, Bakhtin and Kifer \cite{Kifer},
Sowers \cite{Sowers}, Namachchvaya and Sowers~\cite{Namachchvaya-Sowers},
Borodin and Freidlin \cite{Borodin-Freidlin},
Kabanov and Pergamenshchikov \cite{Kabanov-Pergamenshchikov} and references therein.

The specific problem
that we
address in this article is a perturbation of a diffusion in a foliated
manifold $M$ such that the
unperturbed random trajectories lay on the leaves. The perturbations are
taken
transversally to the leaves of the foliation. Here, the slow system is the
transversal component, and the fast system
is given
by the rescaled $y_{{t}/{\varepsilon}}^{\varepsilon}$, where
$y_{t}^{\varepsilon}$
is the solution of the original SDE perturbed by a vector field
$\varepsilon K$.

Our results generalize the recent approach by Li
\cite{Li} on an
averaging principle for a completely
integrable stochastic Hamiltonian system. In this article, as in the classical
approach (see, e.g., \cite{V-Arnold}), Li has
explored the benefits of well-structured geometrical coordinates in the
state
space given by the coordinates of the Liouville torus; these benefits include
vanishing It\^o--Stratonovich correction terms besides also vanishing covariant
derivative of Hamiltonian vector fields in tangent directions to the leaves.
We prove here that an
averaging principle
also holds in a generalized geometrical scope, so that this averaging phenomenon
occurs independently of symplectic
structures (but with possibly slower rates of
convergence).
Compared to Li's
previous result (\cite{Li}, Lemma~3.2), where
the estimates contain a term of order $1/ \sqrt{t}$, our corresponding
estimates
in
Lemma~\ref{lemaequi2} are continuous at $\varepsilon= t=0$.
Some of the rates of convergence of \cite{Li}, Lemma~3.1, are
recovered as particular cases in Corollaries \ref{Corolario_1} and
\ref{Corolario_2}.
In the main result,
we show that in the average, the approximation goes to zero with order
$\varepsilon^\alpha+ \eta( |\ln\varepsilon|^{-{\beta}/{p}})$, where
$\eta$
is the rate of convergence of the unperturbed system to the ergodic
average on
the leaves, $\beta\in(0,1/2)$ and $\alpha\in(0,1)$, so that,
depending on
$\eta$, convergences rate can be much faster or much slower than in
\cite{Li}.

\textit{The set-up.} Let $M$ be a smooth Riemannian manifold with an
$n$-dimen\-sional
smooth foliation; that is, $M$ is endowed with an
integrable regular distribution of dimension $n$; for a definition
and further properties of foliated spaces, see, for example, the
initial chapters of
Tondeur \cite{Tondeur}, Walczak
\cite{Walcak}, among others. We denote by $L_x$ the leaf
of the foliation passing through a point $x\in M$.
For simplicity, we shall assume that the leaves are compact
and that each leaf $L_x$ has a
tubular neighborhood $U\subset M$ where $U$ is diffeomorphic to $L_x
\times
V$, where $V\subset\R^d$ is an open bounded neighborhood of the origin
and $d$
is the codimension of the foliation.
We shall assume an SDE in $M$ whose solution flow preserves the
foliation; that is,
we consider a Stratonovich equation
%
\begin{equation}
\label{eq_original} dx_t = X_0 (x_t) \,dt + \sum
_{k=1}^r X_k
(x_t) \circ dB^k_t,
\end{equation}
where the smooth vector fields $X_k$ are foliated in the
sense that $X_k(x) \in T_x L_x$,
for $k=0,1, \ldots, r$. Here $B_t = (B^1_t, \ldots, B^r_t)$ is a standard
Brownian motion in $\R^r$
with respect to a filtered probability space $(\Omega, \mathcal{F}_t,
\mathcal{F}, \mathbf{P})$. For an
initial condition $x_0$, the trajectories of the solution $x_t$ in this case
lay on the leaf $L_{x_0}$ a.s. Moreover, there exists a (local) stochastic
flow of diffeomorphisms $F_t\dvtx M \rightarrow M$ which, restricted
to the initial leaf, is a flow in the compact $L_{x_0}$.

For a smooth vector field $K$ in $M$, we shall denote the perturbed
system by
$y^{\varepsilon}_t$ which satisfies the SDE
%
\begin{equation}
\label{eq_perturbed} dy^{\varepsilon}_t = X_0
\bigl(y^{\varepsilon}_t\bigr) \,dt + \sum
_{k=1}^r X_k \bigl(y^{\varepsilon}_t
\bigr) \circ dB^k_t + \varepsilon K \bigl(y^{\varepsilon}_t
\bigr) \,dt,
\end{equation}
with the same initial condition $y^{\varepsilon}_0=x_0$ of the unperturbed
system $x_t$.

Our main result, Theorem~\ref{teoremaprincipal}, says that
locally the transversal behavior of
$y^{\varepsilon}_{{t}/{\varepsilon}}$ can be approximated in the
average by an
ordinary differential equation in the transversal space whose
coefficients are given by the average of the transversal component of the
perturbation $K$ with respect to the invariant measure on the leaves
for the
original dynamics of equation (\ref{eq_original}). The reader will
notice by
the end of the proofs that compactness of
the leaves in fact can be substituted by some other boundedness
conditions, added also to some rather technical adjustments which we
will not
address here.
In Sections~\ref{sec2} and \ref{sec3}, we present the main lemmas. The main result
appears in
Section~\ref{sec4}, where we also present a simple illustrative example. In particular,
under some symmetry hypotheses on a foliated system embedded in an Euclidean
space, we use the main theorem to conclude that Lyapunov exponents
in the transversal direction must tend to nonpositive values as
$\varepsilon$
goes to zero;
cf. Proposition~\ref{Proposition Lyapunov
Exp}.

\section{Preliminary results}\label{sec2}

\textit{The coordinate system}.
%
%
Given an initial condition $x_0\in M$, let $U\subset M$ be a
bounded neighborhood of the corresponding leaf
$L_{x_0}$ such that there exists a diffeomorphism $\varphi\dvtx  U
\rightarrow L_{x_0}\times V$, where $V\subset\R^d$ is a connected
open set
containing the origin. The neighborhood $U$ can be considered small enough
such that the closure $\bar{U}\subset M$ and the derivative of
$\varphi
$ is
bounded.
For a fixed diffeomorphism $\varphi$, the space $V$
will be called the \textit{vertical space}. For simplicity, the second (vertical) coordinate of a point
$p \in U$ will be called the vertical projection $\pi(p)\in V$, that is,
$\varphi(p)=(u, \pi(p))$ for some $u\in L_{x_0}$.
Hence for any fixed $v\in V$, the inverse image $\pi^{-1}(v)$ is the compact
leaf $L_x$, where $x$ is any point in $U$ such that the
vertical projection $\pi(x)=v$. In Section~\ref{sec3} below we need the
components of the vertical projection, which shall be denoted by
\[
\pi(p) = \bigl(\pi_1(p), \ldots, \pi_d(p) \bigr) \in V
\subset\R^d
\]
for any $p\in U$.
%
%
Natural examples
of these coordinates systems appear if we consider compact foliation
given by the inverse image of
submersions: values in the image space provide
local coordinates for the vertical space $V$.

%
The next lemma gives information on the order of which the perturbed
trajectories
$y^{\varepsilon}_t$ approach the unperturbed $x_t$ when one varies
$\varepsilon$
and $t$ in
equation (\ref{eq_perturbed}); it will be used to prove that the
dynamics of the rescaled system $y^{\varepsilon}_{ (
{t}/{\varepsilon} )}$ is such
that its time average for any function $g$ in $M$ approximates the time average
of the spacial average of $g$ on the leaves; see Lemma~\ref{lemaequi2}.



%
%

\begin{lemma} \label{lemaequi1} Let $\tau^{\varepsilon}$ be the exit time
of the
process $y^{\varepsilon}_t$ from the neighborhood $U$ of our coordinate system
as above. For any locally Lipschitz continuous
function $f\dvtx M\rightarrow\R$ and $2\leq p < \infty$, we have
\[
\Bigl[\mathbb{E} \Bigl( \sup_{s\leq t \wedge\tau^{\varepsilon}} \bigl\llvert f
\bigl(y^{\varepsilon}_s\bigr)-f(x_s)\bigr\rrvert
^p \Bigr) \Bigr]^{{1}/{p}} \leq K_1 \varepsilon t
e^{K_2 t^p},
\]
where $K_1, K_2\geq0$ are constants depending on upper bounds of the
norms of
the perturbing vector field $K$, on the Lipschitz coefficients of $f$
and on
the derivatives of $X_0, X_1, \ldots,
, X_r$ with respect to the
coordinate system.
\end{lemma}

\begin{pf} Initially write $x_t$ and
$y^{\varepsilon}_t$, the solutions of equations
(\ref{eq_original}) and (\ref{eq_perturbed}), respectively, according
to the
coordinates given by the diffeomorphism $\varphi$: denote
$(u_t, v_t):= \varphi(x_t)$ and
$(u^{\varepsilon}_t,v^{\varepsilon}_t):= \varphi(y^{\varepsilon}_t)$. Then
%
\begin{eqnarray}
\label{primeira des lema1}
\bigl\llvert f\bigl(y^{\varepsilon}_t
\bigr)-f(x_t)\bigr\rrvert &=&\bigl\llvert f\circ
\varphi^{-1} \bigl(u^{\varepsilon}_t,v^{ \varepsilon}_t
\bigr)-f\circ \varphi^{-1} (u_t,v_t)\bigr
\rrvert
\nonumber
\\[-8pt]
\\[-8pt]
\nonumber
&\leq&C \bigl\llvert u^{\varepsilon}_t-u_t\bigr
\rrvert +C\bigl\llvert v^{\varepsilon}_t-v_t\bigr
\rrvert ,
\end{eqnarray}
for some constant $C\geq0 $, where we have used the fact that
$\varphi$ has bounded derivative. The inequality of the statement will be
obtained from inequality (\ref{primeira des lema1}), treating
separately the
norms of the summands above.

For the summand coming from the vertical components, we have that
the unperturbed $v_t \equiv0$.
The perturbing vector field $K$ in our coordinate system is given by $d
\varphi(K)= K_h + K_v$ where $K_h$ and $K_v$ are the horizontal and the
vertical components of $d \varphi(K)$ in the tangent space
$TL_{x_0}\times V$.
The equation for $v^{\varepsilon}_t$ is
given, simply, by
\[
\label{eqcoord2} dv^{\varepsilon}_t= \varepsilon K_v
\bigl(u^{\varepsilon}_t,v^{\varepsilon}_t\bigr)\,dt.
\]
Therefore
%
\begin{eqnarray}
\label{estimativa1}
\sup_{s\leq t \wedge\tau^{\varepsilon}} \bigl\llvert
v^{\varepsilon
}_s - v_s\bigr\rrvert &\leq& \varepsilon\sup
_{s\leq t \wedge\tau^{\varepsilon}} \int^s_0\bigl\llvert
K_v\bigl(u^{\varepsilon}_s,v^{\varepsilon}_s
\bigr)\bigr\rrvert \,ds
\nonumber
\\[-8pt]
\\[-8pt]
\nonumber
&\leq& \varepsilon t \sup_{x\in U} \bigl\llvert K_v(x)
\bigr\rrvert \leq C_1 \varepsilon t,
\end{eqnarray}
where $C_1 = \sup_{x\in U} \llvert K (x)\rrvert $.

For the horizontal summand $\llvert u^{\varepsilon}_t-u_t\rrvert $ in inequality
(\ref{primeira des lema1}), we consider an embedding $i\dvtx M \rightarrow
\R
^N$ of
the compact submanifold $L_{x_0}$
in an Euclidean space $\R^N$ with a
sufficiently large integer $N$.
In $\R^N$ we have, in coordinates, $i(u^{\varepsilon}_t)=(u^{\varepsilon, 1}_t,
\ldots,
u^{\varepsilon,N}_t)$. We look for an equation of $i(u^{\varepsilon}_t) \in
\R^N$.
Let $\tilde{b}_k$ be the vector fields on the embedded
$i(L_{x_0})$ induced by the original vector fields $X_k$ of equation
(\ref{eq_original}), for $k=0, 1, \ldots, r$. Precisely, if $d\varphi
(X_k) =
X_k^h+ X_k^v$ is the decomposition on the horizontal and vertical components,
then $\tilde{b}_k(u,v)= {d} i (X_k^h (u,v) )$ with $u\in i(M)$
and $v\in
V$. By compactness, we extend each $\tilde{b}_k$ to a vector field
$b_k $ in a tubular neighborhood of $i(L_{x_0}) \subset\R^N$ such
that $b_k$
is constant in orthogonal fibres to $i(M)$ in this tubular
neighborhood, for
each $k=0, 1, \ldots, r$. Finally, in canonical coordinates in $\R^N$,
we have,
for
$i=1, \ldots, N$,
%
\begin{equation}
du^{\varepsilon,i}_t = \sum^r_{k=1}b^i_k
\bigl(u^{\varepsilon}_t,v^{\varepsilon
}_t\bigr)\circ
dB^k_t + b^i_0
\bigl(u^{\varepsilon}_t,v^{\varepsilon}_t\bigr) \,dt +
\varepsilon K^i_{h}\bigl(u^{\varepsilon}_t,v^{\varepsilon}_t
\bigr)\,dt. \label{eqcoord1}
\end{equation}

Again, by compactness, the image of the embedding $i(L_{x_0})$ has an induced
metric from $\R^N$ which is uniformly equivalent to the original
metric in
$L_{x_0}$.
Moreover note that, by the choice of the neighborhood $U$, the induced vector
fields $b_0, b_1, \ldots, b_r, K_h$ and their derivatives are bounded. From
equation (\ref{eqcoord1}) we have in each $i$th component, for
$s<\tau^{\varepsilon}$,
%
\begin{eqnarray}
u^{\varepsilon,i}_s - u^i_s &=& \sum
^r_{k=1} \int^s_0{
\bigl(b^i_k\bigl(u^{\varepsilon}_r,v^{\varepsilon
}_r
\bigr)-b^i_k(u_r,v_r)\bigr)\circ
dB^k_r}
\\
&&{} +\int^s_0{\bigl(b^i_0
\bigl(u^{\varepsilon}_r,v^{\varepsilon}_r
\bigr)-b^i_0(u_r,v_r)\bigr)\,dr} +
\varepsilon\int^s_0 {K^i_{h}
\bigl(u^{\varepsilon}_r,v^{\varepsilon}_r\bigr)\,dr}
\label{eq_sem_modulo}.
\end{eqnarray}
In terms of the It\^o integral,
\begin{eqnarray*}
&&\int^s_0 \bigl(b^i_k
\bigl(u^{\varepsilon}_r,v^{\varepsilon}_r
\bigr)-b^i_k(u_r,v_r)\bigr)\circ
dB^k_r \\
&&\qquad =  \int^s_0
\bigl(b^i_k\bigl(u^{\varepsilon}_r,
v^{\varepsilon}_r\bigr)-b^i_k(u_r,v_r)
\bigr)\,dB^k_r
\\
&&\qquad\quad{}+ \frac{1}{2} \int^s_0 \bigl[\nabla
b^i_k \cdot b_k \bigl(u^{\varepsilon}_r,v^{\varepsilon}_r
\bigr) - \nabla b^i_k \cdot b_k
(u_r,v_r) \bigr] \,dr.
\end{eqnarray*}
Hence, taking the absolute values in both sides of equation
(\ref{eq_sem_modulo}) we get, for each
$i$,
%
\begin{eqnarray}
\label{termo_de_correcao_Strat} \bigl\llvert u^{\varepsilon,v}_s - u^i_s
\bigr\rrvert & \leq& \sum^r_{k=1} \biggl
\llvert \int^s_0{\bigl(b^i_k
\bigl(u^{\varepsilon}_r,v^{\varepsilon}_r
\bigr)-b^i_k(u_r,v_r)
\bigr)\,dB^k_r} \biggr\rrvert
\nonumber
\\
&&{} +\frac{1}{2} \sum_{k=1}^r \int
^s_0 \bigl\llvert \nabla b^i_k
\cdot b_k \bigl(u^{\varepsilon}_r,v^{\varepsilon}_r
\bigr) - \nabla b^i_k \cdot b_k
(u_r,v_r) \bigr\rrvert \,dr
\\
&&{} + \int^s_0 \bigl\llvert
b^i_0\bigl(u^{\varepsilon}_r,v^{\varepsilon}_r
\bigr)-b^i_0(u_r,v_r) \bigr
\rrvert \,dr + \varepsilon\int^s_0 \bigl\llvert
K^i_{h}\bigl(u^{\varepsilon}_r,v^{\varepsilon}_r
\bigr) \bigr\rrvert \,dr.\nonumber
\end{eqnarray}
Functions $b^i_0$ and $(\nabla b^i_k \cdot b_k)$ are Lipschitz; hence
for a
common constant $ C_2$,
%
\begin{eqnarray}
\label{mais_zeros} \bigl\llvert u^{\varepsilon,i}_s - u^i_s
\bigr\rrvert &\leq&\Biggl\llvert \sum^r_{k=1}
\int^s_0{\bigl(b^i_k
\bigl(u^{\varepsilon}_r,v^{\varepsilon}_r
\bigr)-b^i_k(u_r,v_r)
\bigr)\,dB^k_r} \Biggr\rrvert
\nonumber
\\[-8pt]
\\[-8pt]
\nonumber
&&{} + C_2 \int^s_0 \bigl\llvert
v^{\varepsilon}_r-v_r\bigr\rrvert \,dr + C_2
\int^s_0 \bigl\llvert u^{\varepsilon}_r
- u_r\bigr\rrvert \,dr + \varepsilon s \sup_{U
}
\llvert K_{h}\rrvert .
\end{eqnarray}
The first deterministic integral, together with inequality
(\ref{estimativa1}) yields
\begin{eqnarray*}
\bigl\llvert u^{\varepsilon,i}_s - u^i_s\bigr
\rrvert & \leq& \sum^r_{k=1} \biggl\llvert
\int^s_0{\bigl(b^i_k
\bigl(u^{\varepsilon}_r,v^{\varepsilon}_r
\bigr)-b^i_k(u_r,v_r)
\bigr)\,dB^k_r} \biggr\rrvert + C_1C_2
\varepsilon s^2
\\
&& {}+ C_2 \int^s_0 \bigl\llvert
u^{\varepsilon}_r - u_r\bigr\rrvert \,dr + C_1
\varepsilon s.
\end{eqnarray*}
Now, for $p\geq1$, there exists a constant $C_3$ such that
%
\begin{eqnarray*}
\bigl\llvert u^{\varepsilon,i}_s - u^i_s
\bigr\rrvert ^p & \leq& C_3 \sum
^r_{k=1} \biggl\llvert \int^s_0{
\bigl(b^i_k\bigl(u^{\varepsilon}_r,v^{\varepsilon}_r
\bigr)-b^i_k(u_r,v_r)
\bigr)\,dB^k_r} \biggr\rrvert ^p +
C_3 \bigl( C_1C_2 \varepsilon s^{2}
\bigr)^p
\\
&&{} + C_3 C^p_2 \biggl( \int
^s_0 \bigl\llvert u^{\varepsilon}_r
- u_r\bigr\rrvert \,dr \biggr)^p + C_3 (
C_1 \varepsilon s )^p.
\end{eqnarray*}
The Cauchy--Schwarz inequality yields
\begin{eqnarray*}
\bigl\llvert u^{\varepsilon,i}_s - u^i_s\bigr
\rrvert ^p & \leq& C_3 \sum^r_{k=1}
\biggl\llvert \int^s_0{\bigl(b^i_k
\bigl(u^{\varepsilon}_r,v^{\varepsilon}_r
\bigr)-b^i_k(u_r,v_r)
\bigr)\,dB^k_r} \biggr\rrvert ^p +
C_3 \bigl( C_1C_2 \varepsilon s^{2}
\bigr)^p
\\
&&{} + C_3 C^p_2 s^{p-1} \int
^s_0 \bigl\llvert u^{\varepsilon}_r
- u_r\bigr\rrvert ^p \,dr + C_3 (
C_1 \varepsilon s )^p.
\end{eqnarray*}
%
Hence
\begin{eqnarray*}
&&\mathbb{E} \sup_{s\leq t \wedge\tau^{\varepsilon}} \bigl\llvert u^{\varepsilon
,i}_s
- u^i_s\bigr\rrvert ^p \\
&&\qquad \leq
C_3 \mathbb{E} \sup_{s\leq t \wedge\tau^{\varepsilon}} \sum
^r_{k=1} \biggl\llvert \int^s_0{
\bigl(b^i_k\bigl(u^{\varepsilon}_r,v^{\varepsilon}_r
\bigr)-b^i_k(u_r,v_r)
\bigr)\,dB^k_r} \biggr\rrvert ^p +
C_3 \bigl( C_1 C_2 \varepsilon
t^{2} \bigr)^p
\\
& &\qquad\quad{} + C_3 C_2^p t^{p-1}
\mathbb{E} \sup_{s\leq t \wedge
\tau^{\varepsilon}} \int^s_0
\bigl\llvert u^{\varepsilon}_r - u_r\bigr\rrvert
^p\,dr + C_3 ( C_1 \varepsilon t )^p
\nonumber
\\
&&\qquad
\leq C_4 \sum^r_{k=1}
\mathbb{E} \biggl[ \int^{t \wedge
\tau^{\varepsilon}}_0{
\bigl(b^i_k\bigl(u^{\varepsilon}_r,v^{\varepsilon
}_r
\bigr)-b^i_k(u_r,v_r)
\bigr)^2 \,dr} \biggr]^{p/2} + C_3 \bigl(
C_1C_2 \varepsilon t^{2} \bigr)^p
\\
&&\qquad\quad{} + C_3 C_2^p t^{p-1} \int
^t_0 \mathbb{E} \Bigl( \sup
_{s\leq
r \wedge\tau^{\varepsilon}} \bigl\llvert u^{\varepsilon}_r -
u_r\bigr\rrvert ^p \Bigr)\,dr + C_3 (
C_1 \varepsilon t )^p,
\end{eqnarray*}
where we have used classical $L^ p$-inequality for martingales (e.g.,
Revuz and
Yor~\cite{Revuz-Yor}).
Using again the Lipchitz property of each
$b_k$ for the terms in the brackets above,
%
\begin{eqnarray}
\label{depois_do_burkholder} && \sum^r_{k=1} \int
^{t \wedge
\tau^{\varepsilon}}_0{\bigl(b^i_k
\bigl(u^{\varepsilon}_r,v^{\varepsilon}_r
\bigr)-b^i_k(u_r,v_r)
\bigr)^2 \,dr}
\nonumber
\\
&&\qquad \leq 2 C^2_2 \biggl( \int^{t \wedge
\tau^{\varepsilon}}_0
\bigl\llvert v^{\varepsilon}_r-v_r\bigr\rrvert
^2\,dr + \int^{t \wedge
\tau^{\varepsilon}}_0\bigl\llvert
u^{\varepsilon}_r-u_r\bigr\rrvert ^2 \,dr
\biggr)
\nonumber
\\[-8pt]
\\[-8pt]
\nonumber
&&\qquad \leq 2 C^2_2 \biggl( \int^t_0
C_1^2 \varepsilon^2 r^2 \,dr + \int
^{t \wedge
\tau^{\varepsilon}}_0 \sup_{s\leq r \wedge\tau^{\varepsilon}} \bigl
\llvert u^{\varepsilon}_r-u_r\bigr\rrvert ^2
\,dr \biggr)
\\
&&\qquad \leq C^2_2 C_1^2
\varepsilon^2 t^3 + 2C^2_2 \int
^{t \wedge
\tau^{\varepsilon}}_0 \sup_{s\leq r \wedge\tau^{\varepsilon}} \bigl
\llvert u^{\varepsilon}_r-u_r\bigr\rrvert ^2
\,dr.\nonumber
\end{eqnarray}
We end up with
%
\begin{eqnarray}
\label{para_p_menor_2} &&\mathbb{E} \sup_{s\leq t \wedge\tau^{\varepsilon}} \bigl\llvert
u^{\varepsilon
,i}_s - u^i_s\bigr\rrvert
^p \nonumber\\
&&\qquad \leq C_4 \sum^r_{k=1}
\mathbb{E} \biggl[ C^2_2 C_1^2
\varepsilon^2 t^3 + 2C^2_2 \int
^{t \wedge
\tau^{\varepsilon}}_0 \sup_{s\leq r \wedge\tau^{\varepsilon}} \bigl
\llvert u^{\varepsilon}_r-u_r\bigr\rrvert ^2
\,dr \biggr]^{p/2}
\nonumber
\\[-8pt]
\\[-8pt]
\nonumber
&&\qquad\quad{} + C_3 C_2^p t^{p-1} \int
^t_0 \mathbb{E} \Bigl( \sup
_{s\leq
r \wedge\tau^{\varepsilon}} \bigl\llvert u^{\varepsilon}_r -
u_r\bigr\rrvert ^p \Bigr)\,dr
\\
&&\qquad\quad{} + C_3 \bigl( C_1C_2 \varepsilon
t^{2} \bigr)^p + C_3 ( C_1
\varepsilon t )^p.\nonumber
\end{eqnarray}
For $p\geq2$ one can use Cauchy--Schwarz again to conclude that there
exists a
constant $C_5$ such that the last expression is less than or equal to
\begin{eqnarray*}
& & C_5 \bigl( C_2 C_1 \varepsilon
t^{3/2} \bigr)^p + C_5 C^p_2
t^{{(p-2)}/{2}} \int^{t \wedge
\tau^{\varepsilon}}_0 \mathbb{E} \sup
_{s\leq r \wedge\tau^{\varepsilon}} \bigl\llvert u^{\varepsilon}_r-u_r
\bigr\rrvert ^p \,dr + C_3 \bigl( C_1C_2
\varepsilon t^{2} \bigr)^p
\nonumber
\\
&&\quad{} + C_3 C_2^p t^{p-1}
\int^t_0 \mathbb{E} \Bigl( \sup
_{s\leq
r \wedge\tau^{\varepsilon}} \bigl\llvert u^{\varepsilon}_r -
u_r\bigr\rrvert ^p \Bigr)\,dr + C_3 (
C_1\varepsilon t )^p
\\
&& \qquad=C_5 \bigl( C_2 C_1 \varepsilon
t^{3/2} \bigr)^p+ C_3 \bigl(
C_1C_2 \varepsilon t^{2} \bigr)^p +
C_3 ( C_1 \varepsilon t )^p
\\
&&\qquad\quad{} + \bigl( C_5 C^p_2 t^{{(p-2)}/{2}}+
C_3 C_2^p t^{p-1} \bigr) \int
^t_0 \mathbb{E} \Bigl( \sup
_{s\leq
r \wedge\tau^{\varepsilon}} \bigl\llvert u^{\varepsilon}_r -
u_r\bigr\rrvert ^p \Bigr)\,dr.
\end{eqnarray*}
Now, summing up over $i$ in the inequalities above leads to
\begin{eqnarray*}
&&\mathbb{E} \sup_{s\leq t \wedge\tau^{\varepsilon}} \bigl\llvert u^{\varepsilon}_r-u_r
\bigr\rrvert ^p \\
&&\qquad \leq C_5 \bigl( C_2
C_1 \varepsilon t^{3/2} \bigr)^p+ C_3
\bigl( C_1C_2 \varepsilon t^{2}
\bigr)^p + C_3 ( C_1 \varepsilon t
)^p
\\
&&\qquad\quad{} + \bigl( C_5 C^p_2 t^{{(p-2)}/{2}}+
C_3 C_2^p t^{p-1} \bigr) \int
^t_0 \mathbb{E} \Bigl( \sup
_{s\leq
r \wedge\tau^{\varepsilon}} \bigl\llvert u^{\varepsilon}_r -
u_r\bigr\rrvert ^p \Bigr)\,dr.
\end{eqnarray*}
We use now the integral form of Gronwall's inequality to find that
\begin{eqnarray*}
\label{desiguap=2} \mathbb{E} \Bigl( \sup_{s\leq t \wedge\tau^{\varepsilon}} \bigl\llvert
u^{\varepsilon}_r-u_r\bigr\rrvert ^p \Bigr)
& \leq& C_6 \varepsilon^p t^p \bigl( 1 +
t^p \bigr) \exp\bigl\{ C_7\bigl(t^{p/2} +
t^p \bigr) \bigr\}
\\
& \leq& C_8 \varepsilon^p t^p \bigl( 1 +
t^p \bigr) \exp\bigl\{ C_9 t^{p} \bigr\}
\nonumber
.
\end{eqnarray*}
Going back to inequality (\ref{primeira des lema1}), now we have
\[
\bigl\llvert f\bigl(y^{\varepsilon}_t\bigr)-f(x_t)
\bigr\rrvert ^p \leq C_{10} \bigl\llvert
v^{\varepsilon}_t-v_t\bigr\rrvert ^p +
C_{10}\bigl\llvert u^{\varepsilon}_t-u_t
\bigr\rrvert ^p.
\]
Hence
\begin{eqnarray*}
\mathbb{E} \Bigl( \sup_{s\leq t \wedge\tau^{\varepsilon}} \bigl\llvert f
\bigl(y^{\varepsilon}_s\bigr)-f(x_s)\bigr\rrvert
^p \Bigr) & \leq& C_{10} \mathbb{E} \sup
_{s\leq t \wedge\tau^{\varepsilon}} \bigl\llvert v^{\varepsilon}_s -
v_s\bigr\rrvert ^p + C_{10} \mathbb{E} \sup
_{s\leq t \wedge\tau^{\varepsilon}} \bigl\llvert u^{\varepsilon}_s-u_s
\bigr\rrvert ^p
\\
&\leq&C_{11} \varepsilon^p t^p +
C_{12}\varepsilon^p t^p\bigl(1+ t^p
\bigr) \exp\bigl(C_9 t^{p}\bigr)
\\
&\leq&C_{13} \varepsilon^p t^p \bigl(1 +
t^p\bigr)\exp\bigl(C_9 t^{p}\bigr).
\end{eqnarray*}
From here, finally, one concludes that there exist constants $K_1$ and $K_2$
such that
\[
\mathbb{E} \Bigl( \sup_{s\leq t \wedge\tau^{\varepsilon}} \bigl\llvert f
\bigl(y^{\varepsilon}_s\bigr)-f(x_s)\bigr\rrvert
^p \Bigr)^{{1}/{p}} \leq K_1 \varepsilon t \exp
\bigl(K_2 t^{p}\bigr).
\]
\upqed\end{pf}

%
%

\textit{On the order of the estimates in the lemma.} Simple
examples show
that the exponential order on the rate of convergence of Lemma~\ref{lemaequi1}
cannot be improved. This fact does not come from the vertical
component, which
is deterministic and easily bounded by $\varepsilon t C_1$; but,
rather, it
comes from the horizontal component. Consider the following example: a
one-dimensional horizontal dynamics which locally is
linear ${d}x_t = x_t \circ{d}B_t$ and a perturbing transversal
vector field $K$ which has constant unitary component on the horizontal
direction, that is, ${d}y^{\varepsilon}_t = y^{\varepsilon}_t \circ
{d}B_t
+ \varepsilon\, {d}t$. The $L^p$ norm of the difference
\[
\E \bigl[x_t - y_t^{\varepsilon} \bigr]^p
\leq\varepsilon\int_0 ^t \E \bigl[ \exp\bigl
\{p(B_t-B_s)\bigr\} \bigr] \,ds
\]
increases exponentially with respect to time $t$.


The next corollary includes the case of a completely integrable
stochastic Hamiltonian system when one uses the action-angle
coordinates; cf.
Li
\cite{Li}, Lemma~3.1.

\begin{corollary} \label{Corolario_1} If the vector fields $X_0,
\ldots
, X_r$
depend only on the
vertical coordinate (null derivative in the directions of the leaves,
as in the Hamiltonian case
\cite{Li}), then the estimates above
can be improved, and for $p\geq1$ there exists a constant $K_1$ such that
\[
\Bigl[\mathbb{E} \Bigl( \sup_{s\leq t \wedge\tau^{\varepsilon}} \bigl\llvert f
\bigl(y^{\varepsilon}_s\bigr)-f(x_s)\bigr\rrvert
^p \Bigr) \Bigr]^{{1}/{p}} \leq K_1 \varepsilon
\bigl(t+ t^ 2\bigr).
\]
\end{corollary}

\begin{pf} In this case the correction term of the Stratonovich stochastic integral
in terms of the It\^o integral in inequality (\ref{termo_de_correcao_Strat})
vanishes, as does the determinist integration of $|u^{\varepsilon}_r -
u_r|$ in inequalities (\ref{mais_zeros}) and (\ref
{depois_do_burkholder}). Hence
inequality (\ref{para_p_menor_2}) improves to
%
\begin{equation}\quad
\mathbb{E} \sup_{s\leq t \wedge\tau^{\varepsilon}} \bigl\llvert u^{\varepsilon
,i}_s
- u^i_s\bigr\rrvert ^p  \leq
C_5 \bigl( C_2 C_1 \varepsilon
t^{3/2} \bigr)^p + C_3 \bigl(
C_1C_2 \varepsilon t^{2} \bigr)^p +
C_3 ( C_1 \varepsilon t )^p.
\end{equation}
The argument in the rest of the proof follows straightforwardly for
$p\geq1$
skipping Gronwall inequality.
\end{pf}

The next corollary includes the case $X_0\equiv0$, cf. \cite{Li}, Lemma~3.1(2),
for stochastic Hamiltonian systems with an action-angle
coordinate system.

\begin{corollary} \label{Corolario_2} If, in addition to the conditions
of Corollary~\ref{Corolario_1} above, we have that the deterministic
vector field $X_0$ is constant when represented with respect to a certain
coordinate system in $U$ (i.e., $b_0$ has null derivative w.r.t. $u$ and~$v$),
then for $p\geq1$ the
estimates can be improved further to $K_1\varepsilon(t+ t^{{3}/{2}})$.
\end{corollary}

\begin{pf} Besides the vanishing terms already mentioned above, the
second
deterministic integral on the right-hand side of inequality
(\ref{termo_de_correcao_Strat}) also vanishes. Hence inequality
(\ref{para_p_menor_2}) simplifies further to
\[
\mathbb{E} \sup_{s\leq t \wedge\tau^{\varepsilon}} \bigl\llvert u^{\varepsilon
,i}_s
- u^i_s\bigr\rrvert ^p  \leq
C_5 \bigl( C_2 C_1 \varepsilon
t^{3/2} \bigr)^p + C_3 ( C_1
\varepsilon t )^p.
\]
\upqed\end{pf}

Yet, from the proof of the Lemma~\ref{lemaequi1}, we have the following:

\begin{remark}\label{remark} For $1\leq p < 2$ and $t$ sufficiently
small, there exist constants $K_1$ and $K_2$ such that
%
\begin{equation}
\label{desig_p} \Bigl[\mathbb{E} \Bigl( \sup_{s\leq t \wedge\tau^{\varepsilon}} \bigl
\llvert f\bigl(y^{\varepsilon}_s\bigr)-f(x_s)\bigr\rrvert
^p \Bigr) \Bigr]^{{1}/{p}} \leq K_1 \varepsilon t
\exp{\bigl( K_2 t^{p}\bigr)}.
\end{equation}
\end{remark}

\begin{pf} One can no longer use Cauchy--Schwarz after inequality
(\ref{para_p_menor_2}). Alternatively, from (\ref{para_p_menor_2}),
use that
\begin{eqnarray*}
&&\mathbb{E} \sup_{s\leq t \wedge\tau^{\varepsilon}} \bigl\llvert u^{\varepsilon
}_s
- u_s\bigr\rrvert ^p \\
&&\qquad \leq C_5 \bigl(
C_2 C_1 \varepsilon t^{3/2} \bigr)^p +
C_5 C^p_2 t^{p/2} \mathbb{E} \sup
_{s\leq t \wedge\tau
^{\varepsilon}} \bigl\llvert u^{\varepsilon
,i}_s -
u^i_s\bigr\rrvert ^p + C_3
\bigl( C_1C_2 \varepsilon t^{2}
\bigr)^p
\\
&&\qquad\quad{} + C_3 C_2^p t^{p-1} \int
^t_0 \mathbb{E} \Bigl( \sup
_{s\leq
r \wedge\tau^{\varepsilon}} \bigl\llvert u^{\varepsilon}_r -
u_r\bigr\rrvert ^p \Bigr)\,dr + C_3
C_1^p \varepsilon^p t^p.
\end{eqnarray*}
If we fix an $0<\delta< 1$ and take $t$ sufficiently small such that
$1-C_5 C^p_2
t^{p/2} > \delta$, then
\begin{eqnarray*}
&&\mathbb{E} \sup_{s\leq t \wedge\tau^{\varepsilon}} \bigl\llvert u^{\varepsilon
}_s
- u_s\bigr\rrvert ^p \\
&&\qquad \leq \delta^{-1}
C_5 \bigl( C_2 C_1 \varepsilon
t^{3/2} \bigr)^p + \delta^{-1} C_3
\bigl( C_1C_2 \varepsilon t^{2}
\bigr)^p
\\
&&\qquad\quad{} + \delta^{-1} C_3 C_2^p
t^{p-1} \int^t_0 \mathbb {E} \Bigl(
\sup_{s\leq
r \wedge\tau^{\varepsilon}} \bigl\llvert u^{\varepsilon}_r -
u_r\bigr\rrvert ^p \Bigr)\,dr + \delta^{-1}
C_3 C_1^p \varepsilon^p
t^p,
\end{eqnarray*}
and the calculation is completed as before using the integral version of
Gronwall's inequality.
\end{pf}

\section{Averaging functions on the leaves}\label{sec3}

Consider a differentiable function $g\dvtx M\rightarrow\mathbb{R}$.
The leaf $L_p$ passing through a point $p\in M$ contains the support of an
invariant measure $\mu_p$ for
the unperturbed system (\ref{eq_original}).
We shall assume that each of such $\mu_p$ is
ergodic. We define a function $Q^g\dvtx V\subset\R^d \rightarrow\R$
such that $Q^g(v)$ is the average of $g$ with
respect to the invariant measure on the corresponding
leaf $\pi^{-1}(v)$. Namely, if $v$ is the vertical
coordinate of $p$, that is,
$\varphi(p)=(u,v)$, then
\[
Q^g (v) = \int_{L_p} g (x) \,d
\mu_p(x).
\]
%
By the ergodic theorem, the function $Q^g$ will also be recovered as
the time
average along the unperturbed trajectories. We are now concerned with
this time
average.

\textit{On the rate of convergence on the leaves.} For deterministic
dynamics it
is known that there are no estimates for the speed of convergence in the
ergodic theorem: it can be as slow as any prescribed
decreasing rate; see, for example, Kakutani and Petersen \cite
{Kakutani-Petersen} and
Krengel \cite{Krengel}. The speed of convergence, in general, depends
both on
the dynamics
and on the function whose average is considered.
In particular, the main theorem in \cite{Krengel} states that given an ergodic
invertible transformation in the interval $[0,1]$, and a prescribed
rate of
convergence (given in terms of a sequence converging to zero which, by
comparison is the prescribed rate of convergence), there exists a
continuous function whose rate of convergence of the ergodic theorem is
slower than the convergence to zero of the prescribed sequence. The result
on the rate of convergence holds almost surely and in the $L^p$ norm.

With this result in mind, one easily construct an example of a
deterministic continuous
dynamics in the torus such that for any prescribed (no matter slow)
rate of
convergence, there exists a function in the torus such that the
rate of convergence of the ergodic theorem is slower than the order of
convergence prescribed.
Indeed, consider the deterministic transformation in the interval
$[0,1]$ given
by a rotation by an irrational angle $\gamma$ in the circle $S^1$
after the
identification of extremal points of the interval [alternatively, $x
\mapsto x+
\gamma (\mathrm{mod}\ 1)$]. Given a prescribed (slow) rate of convergence, Krengel's
theorem guarantees that there exists a continuous function $f(x)$ in $[0,1]$
such that the rate of convergence for this function is slower than the
prescribed speed. Remarks in \cite{Krengel}, page 6, say that $f$ can be
taken continuous in~$S^1$. Represent the flat torus in the square $T^2=
[0,1] \times[0,1]$ with the appropriate identifications. Let
$0=x_0<x_1< x_2<
x_3< x_4= 1$ be a partition of the (horizontal) interval $[0,1]$, whose
subintervals
will be denoted by $I_1$, $I_2$, $I_3$ and~$I_4$. Consider the dynamics in
$T^2$ described in terms of its velocities as: horizontal
component is unitary (constant); it is purely horizontal when
$x\in I_1\cup I_2 \cup I_3$ and makes a rotation in the $y$-coordinate
(vertical) which depends smoothly in $x$, that is, twisting the torus
vertically when $x\in I_4$, such that the total rotation along the interval
$I_4$ is given by the irrational angle $\gamma$.
Consider the continuous function $\tilde{f}(x,y)$ in the torus $T^2$
given by
$\tilde{f}(x,y)= f(y)$, if $x\in I_2 $, $| \tilde{f}(x,y)|$ decay
uniformly in $y$ to zero along $x\in I_3$, $\tilde{f}(x,y)$
vanishes if $x\in I_4$ and $| \tilde{f}(x,y)|$ increase
uniformly in $y$ when $x\in I_1$. Here, by ``increase or decrease
uniformly on
$y$,'' we mean that the function $f(y)$ is multiplied by a factor in
$[0,1]$ which
depends only on $x$. The function $\tilde{f}$ vanishes when the flow performs
the vertical rotation of angle $\gamma$ in the torus. The rate of convergence
of the ergodic
theorem for this continuous deterministic system in $T^2$ is as
slow as the
prescription for the original pure rotation in $S^1$.

The deterministic example above shows that, in our context (which is not
necessarily an elliptic diffusion as in \cite{Li}), there are no optimal
estimates on the
rate of convergence of the ergodic theorem to the averaging function $Q^g$.
We proceed our calculations on the estimates on the rate of
convergence of the averaging principle considering a certain rate of
convergence $\eta(t)$ which will depend on the unperturbed system and
on the
perturbing vector field. We remark that this dependence of an averaging
principle on a certain rate of convergence also had to be used in the
classical
approach as in Freidlin and Wentzell
(\cite{Freidlin-Wentzell}, Chapter~7, Section~9), where they deal with an
averaging principle
with convergence in probability rather than in $L^p$.

The functions which we are interested in averaging on the leaves are the
components of the vertical coordinate of the perturbing vector field
$K$. Hence
function $g$ is going to denote one of the real differentiable functions
$d\pi_1 (K), \ldots, d\pi_d(K)$.

For a fixed initial condition for the unperturbed system $x_0\in U$, denote
the rate of convergence in $L^p$ by $\eta
(x_0, t)$. More precisely, let $\eta(x_0,t)$ be a positive,
asymptotically decreasing to zero real function such that
\[
\biggl[\mathbb{E} \biggl\llvert \frac{1}{ t} \int^{t}_{0}{g
\bigl(F_r (x_0) \bigr)\,dr} - Q^g \bigl(
\pi(x_0 ) \bigr) \biggr\rrvert ^p \biggr]^{{1}/{p}}
\leq\eta(x_0, t),
\]
for all $g= d\pi_i(K)$, with $i=1, \ldots, d$. Continuity of the infinitesimal
generator on the transversal direction and
compactness of the leaves implies that taking supremum, there exists a
rate of
converge $\eta(t)$ in the neighborhood $U$ such that
\[
\eta(x, t) \leq\eta(t),
\]
for any initial condition $x \in U$. Obviously, we have that estimate
$\eta(t)$ is bounded and goes to zero when $t$ goes to infinity.
Examples in
our context include: (1)~uniformly elliptic diffusions as in the Hamiltonian
stochastic systems (cf.
\cite{Li}), where we have $\eta(t) = K/\sqrt{t}$; (2) rotations in
$S^1$ (see
Section~\ref{sec4.1}) $\eta(t)= K/t$; (3) constant averaged function $g$, then
$\eta(t)
\equiv0$.

We assume that $Q^g(v)$ has some degree of continuity with respect to
$v$. This assumption appears in two levels in Lemma~\ref{lema3}: (1) Riemann
integrability of $Q^g (\pi
(y^{\varepsilon}_r))$ with respect to $r$ guarantees the convergence to
zero; (2)
$\alpha$-H\"older continuity guarantees the rate of convergence.

The next lemma estimates the time average of the difference between
function $g$
and $Q^g$ along the
\textit{perturbed} system $y_t^{\varepsilon}$. Here again, the stopping time
$\tau^{\varepsilon}$ denotes the first exit time of the
open neighborhood $U\subset M$ which is diffeomorphic to $L_{x_0}\times
V$. We
have the following estimates for the difference of the averages of
functions $g$ and $Q^g$.

%
%

\begin{lemma}\label{lemaequi2} Let $g$ be one of the functions given
by the
vertical coordinates of the perturbing vector field $K$; that is, let
$g\in
\{d\pi_1(K), \ldots, d\pi_d(K) \}$ and
$Q^g\dvtx V
\rightarrow
\R$ be its corresponding average on the leaves. For $s, t \geq0$ write
\[
\delta(\varepsilon,t) = \int^{(s+t)\wedge\varepsilon\tau^{\varepsilon
}}_{s\wedge
\varepsilon
\tau^{\varepsilon}} g
\bigl(y^{\varepsilon}_{{r}/{\varepsilon}}\bigr) - Q^g\bigl(\pi
\bigl(y^{\varepsilon}_{{r}/{\varepsilon}} \bigr)\bigr) \,dr.
\]
Then $\delta(\varepsilon,t)$ goes to zero when $t$ or $\varepsilon$ tend
to zero.

Moreover, if $Q^g$ is $\alpha$-H\"older continuous, then for
$p\geq
2$, $ 0< \alpha' <\alpha$ and any $\beta\in(0,1/2)$ we have the
following
estimates:
\[
\Bigl(\mathbb{E} \sup_{s\leq
t}\bigl\llvert \delta(\varepsilon,s)
\bigr\rrvert ^p \Bigr)^{{1}/{p}}\leq t \varepsilon^{\alpha'}
h(t,\varepsilon) + t \eta \bigl( t |\ln \varepsilon |^{{(2\beta)}/{p}} \bigr),
\]
where $\eta(t)$ is the rate of convergence in $L^p$ of the ergodic averages
in the unperturbed trajectories on the leaves as defined above, and
$h(t,\varepsilon)$ is continuous for $t, \varepsilon>0$ and converges to
zero when $(t,\varepsilon) \to0$.
\end{lemma}

\begin{pf}
The proof consists of considering a convenient
partition of the interval $(\frac{s}{\varepsilon} \wedge
\tau^{\varepsilon}, \frac{s+t}{\varepsilon} \wedge\tau^{\varepsilon})$
where we
can get the estimates by comparing in each subinterval the average
of the flow of the original system (on the corresponding leaf) with the
average
of the perturbed flow (possibly transversal to the leaves). These
estimates in
each subinterval are obtained using Lemma~\ref{lemaequi1}, so,
a key point in the proof is a careful choice of the increments of such a
convenient partition.

For sufficiently small $\varepsilon$,
we take the following assignment of
increments:
\[
\Delta t = \biggl( \frac{s+t}{|\ln\varepsilon|^{- {2\beta}/{p}}} \biggr)
 \wedge\tau^{\varepsilon} - \biggl(
\frac{s}{|\ln\varepsilon|^{-
{2\beta}/{p}}} \biggr) \wedge\tau^{\varepsilon}.
\]
Hence, the partition $t_n= \frac{s}{\varepsilon} \wedge\tau^{\varepsilon
} + n
\Delta t$, for $0\leq n \leq N$, is such that
\[
\frac{s}{\varepsilon} \wedge\tau^{\varepsilon}=t_0 < t_1 <
\cdots< t_{N} \leq \frac{s+t}{\varepsilon} \wedge\tau^{\varepsilon},
\]
with $N=N(\varepsilon)= \lfloor\varepsilon^{-1} |\ln\varepsilon|^{-{2
\beta}/{p}} \rfloor$
where $\lfloor x\rfloor$ denotes the integer part of $x$.

Initially we represent the left-hand side as the sum
\[
\varepsilon \int^{({(s+t)}/{\varepsilon})\wedge\tau^{\varepsilon}}_{
{s}/{\varepsilon
}\wedge
\tau^{\varepsilon}}{g\bigl(y^{\varepsilon}_r
\bigr)\,dr} = \varepsilon \sum^{N-1}_{n=0}\int
^{t_{n+1}}_{t_n}{g\bigl(y^{\varepsilon}_r
\bigr)\,dr}+ \varepsilon \int^{({(s+t)}/{\varepsilon})\wedge\tau^{\varepsilon
}}_{t_N}{g
\bigl(y^{\varepsilon}_r\bigr)\,dr}.
\]

Denote by $\theta_t$ the canonical shift operator on the probability
space. Let
$ F_t(\cdot,\omega)$ with $ t\geq0$ be the flow of the original unperturbed
system in $M$.
Triangular inequality splits our calculation into four parts,
%
\begin{equation}
\label{split}\bigl |\delta(\varepsilon, t)\bigr|\leq|A_1| + |A_2| +
|A_3|+ |A_4|,
\end{equation}
where
\begin{eqnarray*}
A_1 & = & \varepsilon \sum^{N-1}_{n=0}
\int^{t_{n+1}}_{t_n} \bigl[g\bigl(y^{\varepsilon}_r
\bigr)-g\bigl(F_{r-t_n}\bigl(y^{
\varepsilon}_{t_n},
\theta_{t_n}(\omega)\bigr)\bigr) \bigr]\,dr,
\\
A_2 & = & \varepsilon \sum^{N-1}_{n=0}
\biggl[ \int^{t_{n+1}}_{t_n} g\bigl(F_{r-t_n}
\bigl(y^{
\varepsilon}_{t_n}, \theta_{t_n}(\omega)\bigr)\bigr)
\,dr - \Delta t Q^g \bigl(\pi\bigl(y^{\varepsilon}_{ t_n}
\bigr)\bigr) \biggr],
\\
A_3 & = & \sum^{N-1}_{n=0}
\varepsilon\Delta t Q^g \bigl(\pi\bigl(y^{\varepsilon}_{ t_n}
\bigr)\bigr) - \int^{(s+t)\wedge\varepsilon
\tau^{\varepsilon}}_{s\wedge
\varepsilon
\tau^{\varepsilon}} Q^g
\bigl(\pi\bigl(y^{\varepsilon}_{{r}/{\varepsilon}}\bigr)\bigr)\,dr,
\\
A_4 & = & \varepsilon \int^{({(s+t)}/{\varepsilon})\wedge\tau^{\varepsilon}}_{t_N} g
\bigl(y^{\varepsilon}_r\bigr)\,dr.
\end{eqnarray*}

We proceed by showing that each of the processes $A_1, A_2, A_3$ and
$A_4$ above
tends to zero uniformly on compact intervals. In what follows, we will
explore many times the fact that for $a>0$ and $b\in\R$, $\varepsilon
^a|\ln
\varepsilon|^{b}$ goes to zero as $\varepsilon\searrow0$.
Hence, by construction, except when $\tau^\varepsilon\leq s$ (where,
restricted to
which the lemma is trivial) both $\Delta t$ and
$N$ go to infinity when $\varepsilon$ tends to zero. The increment
$\Delta
t$ is a
random variable bounded by $t |\ln\varepsilon
|^{{2\beta}/{p}}$. Using that the exit time $\tau^\varepsilon$ has a
lower bound
estimate $\tau^\varepsilon\geq\frac{C}{\varepsilon}$ for a positive
constant $C$
independently of $\omega$ [due to equation (\ref{estimativa1})], one
sees that, in fact, the increment $\Delta t = t |\ln\varepsilon
|^{{2\beta}/{p}}$
when $\varepsilon$ is sufficiently small, independently of $\omega$.
%
%

\begin{lemma}\label{lema1} Process $A_1$ converges to zero
uniformly on compact intervals when
$\varepsilon$ goes to zero. More precisely, we have the following
estimates on
the rate of convergence: For any $\gamma\in(0,1)$, there exists a function
$h_1$ such that
\[
\Bigl(\mathbb{E} \sup_{s\leq
t}\llvert A_1 \rrvert
^p \Bigr)^{{1}/{p}}\leq K_1 t \varepsilon
^{\gamma} h_1(t,\varepsilon),
\]
where $h_1$ is continuous in $t, \varepsilon> 0$ and
converges to zero
when
$(t, \varepsilon) \mapsto(0,0)$.
\end{lemma}

\begin{pf} Initially, note that by triangular inequality, and putting the supremum
inside the integral, we get
\begin{eqnarray*}
&&\Bigl(\mathbb{E} \sup_{s\leq
t}\llvert A_1 \rrvert
^p \Bigr)^{{1}/{p}}  \\
&&\qquad\leq \varepsilon \sum
^{N-1}_{n=0} \biggl(\mathbb{E} \biggl[ \int
^{t_{n+1}}_{t_n} \sup_{t_n \leq s\leq
r} \bigl
\llvert g\bigl(y^{\varepsilon}_s\bigr)-g\bigl(F_{s-t_n}
\bigl(y^{
\varepsilon}_{t_n}, \theta_{t_n}(\omega)\bigr)\bigr)
\bigr\rrvert \,dr \biggr]^p \biggr)^{{1}/{p}}.
\end{eqnarray*}
If $ \frac{1}{p}+\frac{1}{q}=1$, by H\"{o}lder's inequality, we have
that the
estimate above is again bounded by
\begin{eqnarray*}
&&\varepsilon\sum^{N-1}_{n=0} \biggl( \mathbb{E}
\biggl[ \biggl( \int^{t_{n+1}}_{t_n} \,dr
\biggr)^{{1}/{q}}\\
&&\hspace*{44pt}{}\times \biggl( \int^{t_{n+1}}_{t_n} \sup
_{t_n \leq s\leq r} \bigl\llvert g\bigl(y^{\varepsilon}_s
\bigr)-g\bigl(F_
{s-t_n}\bigl(y^{\varepsilon}_{t_n},
\theta_{t_n}(\omega)\bigr)\bigr)\bigr\rrvert ^p\,dr
\biggr)^{{1}/{p}} \biggr]^p \biggr)^{{1}/{p}}.
\end{eqnarray*}

In terms of the increments $\Delta t$, the last inequality is again
bounded by
\begin{eqnarray*}
& \leq& \varepsilon(\Delta t)^{{1}/{q}}\sum
^{N-1}_{n=0} \Bigl( \mathbb{E} \Bigl[ \Delta t \sup
_{t_n \leq s\leq t_{n+1}} \bigl\llvert g\bigl(y^{\varepsilon}_s
\bigr)-g\bigl(F_
{s-t_n}\bigl(y^{\varepsilon}_{t_n},
\theta_{t_n}(\omega)\bigr)\bigr)\bigr\rrvert ^p \Bigr]
\Bigr)^{{1}/{p}}
\\
& \leq& \varepsilon \Delta t \sum^{N-1}_{n=0}
\Bigl( \mathbb{E} \Bigl[ \sup_{t_n \leq s\leq t_{n+1}} \bigl\llvert g
\bigl(y^{\varepsilon}_s\bigr)-g\bigl(F_
{s-t_n}
\bigl(y^{\varepsilon}_{t_n}, \theta_{t_n}(\omega)\bigr)\bigr)
\bigr\rrvert ^p \Bigr] \Bigr)^{{1}/{p}}.
\end{eqnarray*}

Lemma~\ref{lemaequi1} (or its corollaries, where the result holds for
$p\geq1$) says that for all $0\leq n \leq N-1$
above, the function $g$ evaluated along trajectories of the perturbed system
compared with $g$ evaluated along the unperturbed trajectories, both
starting at
$y_{t_n}^{\varepsilon}$, satisfies
\[
\Bigl[ \mathbb{E} \sup_{t_n \leq s \leq t_{n+1}} \bigl\llvert g
\bigl(y^{\varepsilon
}_s\bigr) - g(F_{s-t_n}
\bigl(y^{\varepsilon}_{t_n},\theta_{
t_n}(\omega)\bigr) \bigr
\rrvert ^p \Bigr]^{{1}/{p}} \leq K_1 \varepsilon \Delta
t e^{K_2 (\Delta
t)^{p}}.
\]
Hence, using that $N$
has order $[\varepsilon^{-1} |\ln\varepsilon
|^{-{2\beta}/{p}}]$ and that $\Delta t \leq t |\ln\varepsilon
|^{{2\beta}/{p}}$,
\begin{eqnarray*}
\Bigl[\mathbb{E} \sup_{s\leq
t} |A_1|^p
\Bigr]^{{1}/{p}}& \leq& K_1 N \varepsilon^2 (\Delta
t)^2 e^{K_2 (\Delta t)^{p}}
\\
&\leq& K_1 \varepsilon^2 \bigl[ \varepsilon^{-1} |
\ln \varepsilon|^{-{2\beta}/{p}}\bigr] t^2 |\ln\varepsilon|^{{4\beta}/{p}}
e^{K_2 (t |\ln
\varepsilon|^{{2 \beta}/{p}})^{p}}
\\
&= & K_1 t \varepsilon^{\gamma} h_1 (
\varepsilon, t)
\end{eqnarray*}
for any $\gamma\in(0,1) $ where
\[
h_1 (\varepsilon, t) = t \varepsilon^{{(1- \gamma)}/{2}}  |\ln
\varepsilon|^{{2\beta}/{p}} \exp \biggl\{ \biggl(\frac{1- \gamma}{2} \biggr) \ln
\varepsilon+ K_2 t^{p} |\ln\varepsilon|^{2 \beta} \biggr\},
\]
which satisfies the required properties for $\beta\in(0, 1/2)$.
\end{pf}

%
%

%
%

%
%

\begin{lemma}\label{lema2} Process $A_2$ in equation
(\ref{split}) goes to zero with the following rate of convergence:
\[
\Bigl[\mathbb{E} \sup_{s\leq
t} |A_2|^p
\Bigr]^{{1}/{p}} \leq t \eta \bigl( t |\ln\varepsilon |^{{2\beta}/{p}} \bigr),
\]
where $\eta(t)$ is the rate of convergence in $L^p$ of the ergodic averages
in the unperturbed trajectories on the leaves.
\end{lemma}

\begin{pf} We have
\begin{eqnarray*}
&&\Bigl[\mathbb{E} \sup_{s\leq
t} |A_2|^p
\Bigr]^{{1}/{p}}\\
&&\qquad \leq \varepsilon \Biggl[\mathbb{E} \Biggl\llvert \sum
^{N-1}_{n=0}\biggl[ \int^{t_{n+1}}_{t_n}{g
\bigl(F_{r-t_n}\bigl(y^{\varepsilon}_{t_n},
\theta_{t_n } (\omega)\bigr)\bigr)\,dr}- \Delta t Q^g \bigl(
\pi\bigl(y^{\varepsilon}_{ t_n}\bigr)\bigr) \biggr] \Biggr\rrvert
^p \Biggr]^{{1}/{p}}
\\
&&\qquad \leq \varepsilon \sum^{N-1}_{n=0}
\biggl[\mathbb{E} \biggl\llvert \int^{t_{n+1}}_{t_n}{g
\bigl(F_{r-t_n}\bigl(y^{\varepsilon}_{t_n},
\theta_{t_n } (\omega)\bigr)\bigr)\,dr}- \Delta t Q^g \bigl(
\pi\bigl(y^{\varepsilon}_{ t_n}\bigr)\bigr) \biggr\rrvert ^p
\biggr]^{{1}/{p}}
\\
&&\qquad =  \varepsilon\Delta t \sum^{N-1}_{n=0}
\biggl[\mathbb{E} \biggl\llvert \frac{1}{\Delta t} \int^{t_{n+1}}_{t_n}{g
\bigl(F_{r-t_n}\bigl(y^{\varepsilon}_{t_n},
\theta_{t_n } (\omega)\bigr)\bigr)\,dr}- Q^g \bigl(\pi
\bigl(y^{\varepsilon}_{ t_n}\bigr)\bigr) \biggr\rrvert ^p
\biggr]^{{1}/{p}}.
\end{eqnarray*}
For all $n=0, \ldots, N-1$, by the ergodic theorem, the two
terms inside the
modulus converge to each other when $\Delta t$ goes to infinity with rate
$\eta(\Delta t)$. Hence, for small~$\varepsilon$,
we have
\begin{eqnarray*}
\Bigl[\mathbb{E} \sup_{s\leq
t} |A_2|^p
\Bigr]^{{1}/{p}}& \leq& \varepsilon N (\Delta t) \eta(\Delta t)
\\
& \leq& \varepsilon \bigl[ \varepsilon^{-1} |\ln \varepsilon|^{-{2\beta}/{p}}
\bigr] t |\ln\varepsilon|^{{2\beta}/{p}} \eta \bigl( t |\ln\varepsilon |^{{2\beta}/{p}}
\bigr)
\\
& =& t \eta \bigl( t |\ln\varepsilon|^{{2\beta}/{p}} \bigr).
\end{eqnarray*}
\upqed\end{pf}

%
%

\begin{lemma}\label{lema3} $A_3$
converges to zero when $t$ or $\varepsilon$ go to $0$. Moreover, if $Q^g$ is
$\alpha$-H\"older continuous
then the rate
of
convergence satisfies
\[
\Bigl(\mathbb{E} \sup_{s\leq
t}\llvert A_3 \rrvert
^p \Bigr)^{{1}/{p}}\leq C \varepsilon^{\alpha}
t^{1+ \alpha} |\ln \varepsilon|^{{2\beta\alpha}/{p}},
\]
for a positive constant $C$.
\end{lemma}

\begin{pf} Consider the partition given by the sequence $(\varepsilon
t_n)_{0\leq
n \leq N}$ of the interval $( s \wedge
\varepsilon\tau^{\varepsilon}, (s+ t) \wedge
\varepsilon\tau^{\varepsilon})$.
The convergence to zero here corresponds to the existence of
the Riemann integral.
Moreover, assuming that $Q^g$ is
$\alpha$-H\"older continuous, then
%
%
\begin{eqnarray}\label{Ineq:
proof lemma
3}
|A_3| &\leq& \sum^{N-1}_{n=0}
\varepsilon\Delta t \sup_{\varepsilon t_n < r
\leq
\varepsilon t_{n+1}}\bigl|Q^g\bigl(\pi
\bigl(y^{\varepsilon}_{ t_n}\bigr)\bigr)- Q^g\bigl(\pi
\bigl(y^{\varepsilon}_{{r}/{\varepsilon}}\bigr)\bigr) \bigr|
\nonumber
\\[-8pt]
\\[-8pt]
\nonumber
&\leq& \varepsilon\Delta t \sum^{N-1}_{n=0}
C_1 \Bigl( \sup_{\varepsilon
t_n < r
\leq
\varepsilon t_{n+1}}\bigl| \pi\bigl(y^{\varepsilon}_{ {\varepsilon t_n}/{\varepsilon}}
\bigr)- \pi \bigl(y^{\varepsilon}_{{r}/{\varepsilon}}\bigr) \bigr| \Bigr)^{\alpha}
.
\end{eqnarray}
We use now that
\[
\bigl| \pi\bigl(y^{\varepsilon}_{ {s}/{\varepsilon}}\bigr)- \pi
\bigl(y^{\varepsilon}_{{t}/{\varepsilon}}
\bigr)\bigr | \leq \Bigl( \sup_{x\in U} \bigl|K(x)\bigr| \Bigr) |s-t|
\]
for all $s,t \geq0$, such that the right-hand side is independent of
$\varepsilon$
[we are going to use this fact again, which follows by equation
(\ref{estimativa1}); for more details, see the very beginning of
proof of Theorem~\ref{teoremaprincipal}]. Hence, continuing the estimates for $|A_3|$,
inequality
(\ref{Ineq: proof lemma 3}) above implies that
\begin{eqnarray*}
|A_3| &\leq& C_2 N (\varepsilon\Delta t)^{(1+ \alpha)}
\\
& \leq& C_2 \bigl( \varepsilon^{-1} |\ln
\varepsilon|^{-{2 \beta}/{p}} \bigr) \bigl( \varepsilon t |\ln \varepsilon|
^{{2\beta}/{p}}
\bigr) ^{1 + \alpha}
\\
& \leq& C_2 \varepsilon^{\alpha} t^{1 + \alpha} |\ln
\varepsilon|^{{2\beta\alpha}/{p}}
\end{eqnarray*}
for a constant $C_2$.
\end{pf}

%
%

\begin{lemma}\label{lema4} Process $A_4$ converges to zero
with
\[
\Bigl(\mathbb{E} \sup_{s\leq
t}\llvert A_4 \rrvert
^p \Bigr)^{{1}/{p}} \leq C t \varepsilon |\ln \varepsilon|^{{2\beta}/{p}}.
\]
\end{lemma}

\begin{pf} Denote
\[
C = \sup_{x\in U} \bigl|g(x)\bigr|.
\]
The result follows straightforwardly since
\[
\varepsilon\biggl\llvert \int^{({(s+t)}/{\varepsilon})\wedge
\tau^{\varepsilon}}_{t_N}{g
\bigl(y^{\varepsilon}_r\bigr)\,dr}\biggr\rrvert \leq C\varepsilon\Delta t
= C t \varepsilon|\ln\varepsilon|^{{2\beta}/{p}}.
\]
\upqed\end{pf}

Now, going back to the proof of Lemma~\ref{lemaequi2}, the statement
follows by inequality~(\ref{split}) and adding up the estimates of the last
four lemmas (\ref{lema1} up to \ref{lema4}).
\end{pf}

%
%

\section{An averaging principle}\label{sec4}

We state the averaging principle in the next theorem. To use Lemma~\ref{lemaequi2} of the previous section we have to assume
regularity in the average function $Q^g$, which naturally depends on
$g$, on the
foliated coordinate system and on the
transversal behavior of the invariant measures on the leaves of the original
foliated system. We are going to assume the following condition:

\renewcommand{\thehypo}{(H)}
\begin{hypo}\label{hypoH}
Let $g$ be one of the functions given
by the
vertical coordinates of the perturbing vector field $K$, that is, $g\in
\{d\pi_1(K), \ldots, d\pi_d(K) \}$. We assume that its corresponding
average on
the leaves $Q^g\dvtx V
\rightarrow
\R$ is Lipschitz.
\end{hypo}

This hypothesis holds naturally if the invariant measures $\mu_p$ for the
unperturbed foliated system has a sort of weakly continuity on $p$. For
deterministic
systems it corresponds to a certain regularity in the sense that there
is no
bifurcation
with respect to the vertical parameter $v\in V$.

Let $v(t)$ be the solution of the deterministic ODE in the transversal
component $V\subset\R^n$
%
\begin{equation}
\frac{dv}{dt}= \bigl(Q^{ {d\pi_1 (K)}}, \ldots, Q^{d\pi_d(K)}\bigr)
\bigl( v(t) \bigr)
\end{equation}
with initial condition $v(0)=\pi(x_0)=0$. Let $T_0$ be the time that
$v(t)$
reaches the boundary of $V$.

%
%

\begin{theorem} \label{teoremaprincipal} Assuming Hypothesis~\ref{hypoH} above,
we have:
\begin{longlist}[(1)]
\item[(1)] For any $0<t<T_0$, $\beta\in(0,1/2)$, $\alpha\in(0,1)$
and $2
\leq p< \infty$ ($1 \leq p$ if in the conditions of Corollaries
\ref{Corolario_1}, \ref{Corolario_2} or Remark~\ref{remark}), there
exist functions $C_1=C_1(t)$ and $C_2=C_2(t)$
such that
\[
\Bigl[\mathbb{E} \Bigl( \sup_{s\leq t}\bigl\llvert \pi \bigl(
y^{\varepsilon}_{ ( {s}/{\varepsilon} ) \wedge\tau
^{\varepsilon}} \bigr)-v (s)\bigr\rrvert ^p \Bigr)
\Bigr]^{{1}/{p}} \leq C_ 1 \varepsilon^{\alpha} +
C_2 \eta \bigl( t |\ln \varepsilon|^{{2\beta}/{p}} \bigr),
\]
where $\eta(t)$ is the rate of convergence in $L^p$ of the
ergodic averages
of the unperturbed trajectories on the leaves.

\item[(2)] For $\gamma> 0 $, let
\[
T_{\gamma} = \inf \bigl\{t>0 |\operatorname{dist} \bigl(v(t), \partial V \bigr)
\leq \gamma\bigr\}.
\]

The exit times of the two systems satisfy the estimates
\[
\mathbb{P}\bigl(\varepsilon\tau^{\varepsilon}<T_{\gamma}\bigr)\leq
\gamma^{-p} \bigl[ C_ 1(T_{\gamma})
\varepsilon^{\alpha} + C_2 (T_{\gamma}) \eta \bigl(
T_\gamma|\ln \varepsilon|^{{2\beta}/{p}} \bigr) \bigr]^p.
\]
\end{longlist}
\end{theorem}

The second part of the theorem above guarantees the robustness of the averaging
phenomenon
in the transversal direction.

\begin{pf*}{Proof of Theorem \ref{teoremaprincipal}} The gradient of each real function $\pi_i$
is orthogonal to the leaves, hence by It\^o's formula, for
$i=1,2,\ldots,
d$, we have that
\[
\pi_i \bigl( y^{\varepsilon}_{({t}/{\varepsilon})\wedge\tau^{\varepsilon}} \bigr) = \int
^{t\wedge
\varepsilon\tau^{\varepsilon}}_0 d\pi_i(K)
\bigl(y^{\varepsilon}_{
{s}/{\varepsilon
}}\bigr) \,ds.
\]
Lemma~\ref{lemaequi2} for the function $d \pi_i(K)$ in $M$,
triangular inequality and Hypothesis~\ref{hypoH} imply that
\begin{eqnarray*}
&&\bigl\llvert \pi_i \bigl( y^{\varepsilon}_{({t}/{\varepsilon})\wedge\tau^{\varepsilon}}
\bigr)-v_i (t)\bigr\rrvert \\
&&\qquad\leq \int^{t\wedge
\varepsilon\tau^{\varepsilon}}_0{
\bigl\llvert Q^{d\pi_i(K)} \bigl( \pi \bigl( y^{\varepsilon}_{({s}/{\varepsilon})\wedge
\tau^{\varepsilon}}
\bigr) \bigr)-Q^{d\pi_i(K)} \bigl( v(s)\bigr)\bigr\rrvert \,ds} +\bigl|
\delta_i(\varepsilon,t)\bigr|
\\
&&\qquad\leq C_i \int^t_0
\bigl\llvert \pi \bigl( y^{\varepsilon}_{({s}/{\varepsilon})\wedge\tau^{\varepsilon}} \bigr)-v (s) \bigr\rrvert
\,ds +\bigl| \delta_i( \varepsilon,t)\bigr|,
\end{eqnarray*}
where each $C_i$ is the Lipschitz constant of $Q^{d\pi_i(K)}$, and
$\delta_i(\varepsilon,
t)$ is defined in Lemma~\ref{lemaequi2}.
Summing up the $i$'s and using Gronwall's lemma, we have, for a
constant $C$,
\[
\bigl\llvert \pi \bigl( y^{\varepsilon}_{({t}/{\varepsilon})\wedge\tau^{\varepsilon}} \bigr)-v (t)\bigr
\rrvert \leq e^{Ct} \sum_{i=1}^n
\bigl|\delta_i( \varepsilon,t)\bigr|.
\]
The first part of the theorem follows by Lemma~\ref{lemaequi2}.

For the second part we have the following estimates:
\begin{eqnarray*}
\mathbb{P}\bigl(\varepsilon\tau^{\varepsilon}<T_{\gamma}\bigr) &\leq&
\mathbb{P} \Bigl( \sup_{s\leq
T_{\gamma}\wedge
\varepsilon\tau^{\varepsilon}}\bigl\llvert v(s)-\pi \bigl(
y^{\varepsilon}_{({s}/{\varepsilon})\wedge\tau^{\varepsilon}} \bigr)\bigr\rrvert >\gamma \Bigr)
\\
&\leq& \gamma^{-p}\mathbb{E} \Bigl( \sup
_{s\leq
T_{\gamma}\wedge
\varepsilon\tau^{\varepsilon}} \bigl\llvert v(s)-\pi \bigl( y^{\varepsilon}_{({s}/{\varepsilon})\wedge\tau^{\varepsilon}}
\bigr)\bigr\rrvert ^p \Bigr)
\\
&\leq& \gamma^{-p} \bigl[ C_ 1(T_{\gamma})
\varepsilon^{\alpha} + C_2 (T_{\gamma}) \eta \bigl(
T_\gamma|\ln \varepsilon|^{{2\beta}/{p}} \bigr) \bigr]^p.
\end{eqnarray*}
\upqed\end{pf*}

\subsection{A detailed example}\label{sec4.1}

The following simple example illustrates the
framework
where the averaging principle described in this section holds. Consider
$M= \R^3
- \{ (0,0,z), z\in R \}$
with the 1-dimension horizontal circular foliation of $M$ where the leaf
passing through a point $p=(x,y,z)$ is given by the circle $L_p = \{ (
\sqrt{x^2+ y^2}
\cos\theta, \sqrt{x^2+ y^2} \sin
\theta, z)$, $\theta\in[0,2 \pi] \}$. Consider the foliated linear
SDE on
$M$ consisting of random rotations
\[
dx_t = \pmatrix{ 0 & -1 & 0
\cr
1 & 0 & 0
\cr
0 & 0 & 0 } x_t (
\lambda_ 1 \,dt + \lambda_ 2 \circ dB_t).
\]

For an initial condition $p_0=(x_0,y_0,z_0)$, say with $x_0 \geq0$, consider
the local foliated coordinates in the neighborhood $U= \R^3
\setminus\{(x, 0, z); x\leq0; z\in\R\}$ given by cylindrical coordinates.
Hence, using
the same notation as before, $\varphi= (u,v)$ will be defined by
$\varphi\dvtx U
\subset
M \rightarrow(-\pi, \pi) \times\R_{>0}\times\R$, where $u \in
(-\pi
, \pi)$
is angular and $v=(r,z) \in
\R_{>0}\times\R$ such that $\varphi^{-1}\dvtx (u, v)\mapsto(r \cos u, r
\sin u, z)
\in M$. In this coordinates system, the transversal projections $\pi
_1$ and
$\pi_2$
correspond to the radial $r$-component and the $z$-coordinate,
respectively.

For $\lambda_1, \lambda_2\in\R$ with $|\lambda_1|+ |\lambda_2|>
0$, the
invariant measures $\mu_p$ in the leaves $L_p$ passing through points
$p\in M$
are
given by normalized Lebesgue measures in $L_p$, which here corresponds
to the normalized angle 1-form. Note that Hypothesis~\ref{hypoH} is satisfied.
We
investigate the following
effective behavior of a small
transversal perturbation of order $\varepsilon$:
\[
dy^{\varepsilon}_t = \pmatrix{ 0 & -1
& 0
\cr
1 & 0 & 0
\cr
0 & 0 & 0 } y^{\varepsilon}_t (
\lambda_ 1 \,dt + \lambda_ 2 \circ dB_t) +
\varepsilon K \bigl(y^{\varepsilon}_t\bigr)\,dt,
\]
with initial condition $x_0=
(1,0,0)$.
Typically, due to uniform ellipticity, for any perturbing vector field, the
average of its transversal component on the leaves converges with rates at
least $
\frac{1}{\sqrt{t}}$. Further, if $\lambda_1\neq0$
and $\lambda_2 = 0$, due to periodicity, the rate of convergence improves
to order $\eta(t)\sim\frac{1}{t}$.
Here we consider two classes of
perturbing vector
field
$K$.

\textit{Constant perturbation}. (A) Assume that the perturbation is given
by
a vector
field which is constant $K=(k_1, k_2,
k_3)$ with respect to Euclidean coordinates in $M$. In this case we
have that the convergence of the time average on the leaves to the
ergodic average is given by $\eta(t)\equiv0$. Initially, to fix the
ideas, assume that $k_3=0$. Then, not
only the average on the $z$-component vanishes,
that is, $Q^{d\pi_2 K }=0$, but also, by the geometrical symmetry of
$K$ with
respect to the invariant measure, the average radial $r$-component
also vanishes, that is, $Q^{d\pi_1 K}=0$. Hence the transversal
component
in the main theorem is constant
$v(t)=(r(0), z(0))$ for all $t\geq0$.

The vertical components in the statement of the first formula in Theorem~\ref{teoremaprincipal}, in this example, concern only the difference
between the initial radius $r(0)=1$ and
$r(\frac{t}{\varepsilon}\wedge\tau^{\varepsilon})= \pi_1 (
y^{\varepsilon}_{({t}/{\varepsilon})\wedge\tau^{\varepsilon}}
)$. Precisely, we have that
\[
\biggl[\mathbb{E} \biggl( \sup_{s\leq t} \biggl\llvert r \biggl(
\frac{s}{\varepsilon}\wedge \tau^{\varepsilon} \biggr) - 1 \biggr\rrvert
^p \biggr) \biggr]^{{1}/{p}}
\]
goes to zero, for a fixed $t$, with order $
\varepsilon^\alpha$, with $\alpha\in(0,1)$.

(B) Vertical perturbations. Now, for constant and vertical $K = (0,0,
k_3)$, the radial average $Q^{d \pi_1 K}$ is null, but $Q^{d \pi_2
K}$ equals $k_3$ for every leaf in $M$. Hence the averaged system
$v(t)= (r(0),
k_3t)$ is constant in the radial component and increases linearly in the
$z$-coordinate. The perturbed
systems have the simple solution
\[
y^{ \varepsilon}_{{t}/{\varepsilon}} = \pmatrix{
\displaystyle\cos \biggl( \frac{\lambda_ 1 t}{\varepsilon} + \lambda_2 B_{{t}/{\varepsilon}} \biggr)
\cr
\displaystyle\sin \biggl( \frac{\lambda_ 1 t}{\varepsilon} + \lambda_2
B_{{t}/{\varepsilon}} \biggr)
\cr
k_3 t }.
\]
Hence, the comparison
\[
\bigl| \pi_2 \bigl( y^{\varepsilon}_{({t}/{\varepsilon})\wedge\tau^{\varepsilon}} \bigr) - v(t)\bigr|
\equiv0
\]
for all $t\geq0$ and the convergence to zero is trivial.

\textit{Linear perturbation}.
Consider a linear perturbation of the
form $K(x,y,z)= (x, 0,0 )$. In this case, again, the $z$-coordinate average
vanishes trivially. For the radial component, we have that $d \pi_1 K
= r_0
\cos^2 u$,
where $u$ is the angular coordinate of $p$ whose distance to the $z$-axis
($r$-coordinate) is $r_0$. Hence the average with respect to the invariant
measure on the leaves is given by $ Q^{d\pi_1 K}= r/2$ for leaves with radius
$r$.
%
The transversal system stated in the theorem is then $v(t)= (e^
{{t}/{2}} r(0),
z(0))$. Hence the result guarantees that the radial part of
$y^{\varepsilon}_{({t}/{\varepsilon})\wedge\tau^{\varepsilon}}$ must
have a behavior
close to the exponential $e^{{t}/{2}}$ in the sense that
\[
\biggl[\mathbb{E} \biggl( \sup_{s\leq t} \biggl\llvert r \biggl(
\frac{t}{\varepsilon} \biggr) - e^{
{t}/{2}}\biggr\rrvert ^p \biggr)
\biggr]^{{1}/{p}}
\]
goes to zero when $\varepsilon$ goes to zero.

The fundamental
solution of the linear perturbed Stratonovich systems $y_t^{\varepsilon}$
is given
by the
exponential
of the matrix for each fixed $t$
\[
\pmatrix{ \varepsilon t & - (
\lambda_1 t+ \lambda_2 B_t) & 0
\cr
\lambda_1 t+ \lambda_2 B_t & 0 & 0
\cr
0 & 0 & 0 },
\]
where the eigenvalues corresponding to the first two coordinates (horizontal
plane)
are
\[
\lambda_{1,2} = \frac{\varepsilon t \pm\sqrt{ \varepsilon^2 t^2 - 4
(\lambda
_1 t+
\lambda_2 B_t)^2}}{2},
\]
whose real part is given by $\varepsilon t/2$ with probability
increasing to 1 as $\varepsilon$ goes to
zero. This exponential rates coincides with that one above, guaranteed
by the
Theorem~\ref{teoremaprincipal}.

\textit{Lyapunov exponents}. For foliated manifolds embedded in $\R
^N$, the
symmetry of the perturbing vector
field $K$ with respect to the geometry of the leaves, hence also with
respect to Lebesgue invariant measure, as presented in the case of
constant~$K$,
has implied that the transversal average $Q^{d\pi K}$
vanishes. This phenomenon also appears in a couple of other examples
where the
leaves are not only diffeomorphic to each other, but also have this
symmetry in
the
sense that the integration of a constant perturbation $K\in\R^N$ with respect
to the Lebesgue measure is zero. Here we mention a
couple of
simple examples: the spherical foliation of $\R^n \setminus\{0\}$, nested
torus (increasing the smaller radius)
foliation of the solid torus minus the central circle $S^1 \times D^ 2
\setminus S^{1} \subset\R^3$ or more generally (when they exists) tubular
foliation of
$\R^n \setminus\{ C \}$, with $C$ a compact set (this context also includes
the
Hamiltonian case with the Lyouville foliation of the symplectic
structure of
$\R^{2n}$ as in \cite{Li}). In these symmetric geometrical
configuration, if
the invariant measure on the leaves are the Lebesgue measures (taking gradient
Brownian motion on the leaves, e.g., as in~\cite{CRL-fgBm}) the averages
of
$Q^{d\pi}$ vanish. Hence our main theorem says that on the average, the
trajectories of the perturbed system stay somehow close to the initial
leaf as
$\varepsilon$ decreases to zero.

The Lyapunov exponent of the system in the direction of a tangent
vector $v\in
T_{x_0}M$ contains the long time behavior of points close to $x_0$ in the
direction of $v$; for details on the definition, properties,
existence
conditions, multiplicative ergodic theory, etc., see, for example,
Arnold \cite{Ludwig} and the references therein, among many others. In
particular, under
the symmetric geometrical circumstances above, if there exists the Lyapunov
exponents of the perturbed system $y_t^{\varepsilon}$, Theorem~\ref{teoremaprincipal} will imply that in transversal directions the Lyapunov
exponent cannot be too far from zero. This vanishing property must
happen with
multiplicity given by the codimension of the foliation, as in the
examples of
the paragraphs above, where the asymptotic relevant
parameters (rotation
number \cite{Arn-Imkeller,CRL-RN,Ruffino-SD}, and Lyapunov
exponents) do exist. In short:

\begin{proposition}[(Continuity of Lyapunov exponents)] \label{Proposition
Lyapunov
Exp} Assume that the perturbed system $y^{\varepsilon}_t$ does have
Lyapunov exponents a.s. at the assigned initial condition. If the averaged
perturbation in the leaves vanishes, that is, $Q^{d \pi K}=0$, then a
number, given
by the codimension of the foliation, of Lyapunov exponents in the spectrum
go to nonpositive values as $\varepsilon$ goes to zero.
\end{proposition}

\begin{pf} In fact, Theorem~\ref{teoremaprincipal} says that the perturbed system\break
increases~in~the transversal coordinates with order $| \pi(y^{\varepsilon}_t)|
\lesssim C_1 (\varepsilon t) \varepsilon^\alpha+\break C_2 (\varepsilon t)
\eta ( (\varepsilon t) |\ln\varepsilon|^{{2\beta}/{p}} )$.
Hence, any
exponential
behavior, if it exists, tends to nonpositive values as $\varepsilon$
decreases to
zero.
\end{pf}

\section*{Acknowledgments} This article was written while the
second author's was in a sabbatical visit to University of Humboldt. He would
like to express his gratitude to Professor Peter Imkeller and his
research group for
their nice and friendly hospitality.

\printaddresses
\end{document}